\documentclass[reqno,11pt]{amsart}
\usepackage{amsmath,amsthm,amssymb,amsfonts,amscd}
\usepackage{graphicx,color}

\input xy
\xyoption{all}

\newcommand{\R}{{\mathbb{R}}}
\newcommand{\Z}{{\mathbb{Z}}}
\newcommand{\C}{{\mathbb{C}}}

\def \cA{{\mathcal A}}
\def \cS{{\mathcal S}}

\def \delm#1{{\delta_{[#1]}}}

\def \0{{\bf 0}}
\def \1{{\bf 1}} 
\def \2i{\frac{\sqrt{-1}}{2}}

\def \Mat{\mathit{Mat}}

\def \nabb{{\bar \nabla}}

\def \Th{{\hat T}}

\def \yh{{\hat y}}

\def \fe{{\frak e}}

\newcommand{\ti}{\tilde}
\newcommand{\wti}{\widetilde}

\def \lpartial#1{\overleftarrow{\partial_{#1}}}
\def \rpartial#1{\overrightarrow{\partial_{#1}}}

\def \fpartial#1{\frac{\partial}{\partial {#1}}}

\def \flpartial#1{\frac{\overleftarrow{\partial}}{\partial #1}}
\def \frpartial#1{\frac{\overrightarrow{\partial}}{\partial #1}}

\newtheorem{thm}{Theorem}[section]

\newtheorem{lem}[thm]{Lemma}

\theoremstyle{definition}

\newtheorem{defn}[thm]{Definition}

\newtheorem*{pf}{Proof}

\numberwithin{equation}{section}
\newcommand{\bp}{\begin{pmatrix}}
\newcommand{\ep}{\end{pmatrix}}
\newcommand{\bps}{\begin{smallmatrix}}
\newcommand{\eps}{\end{smallmatrix}}

\def\Hom{\mathrm{Hom}}

\def \Mat{\mathrm{Mat}}
\def \Ker{\mathrm{Ker}}

\def \Ob{\mathrm{Ob}}

\def \ov#1{\frac{1}{#1}}

\begin{document}
\hfill YITP-05-64\\

\title{Star product formula of theta functions}

\date{October, 2005}
\author{Hiroshige Kajiura}
\address{Yukawa Institute for Theoretical Physcics, Kyoto University,
Kyoto 606-8502, Japan}
\email{kajiura@yukawa.kyoto-u.ac.jp}

\begin{abstract}

As a noncommutative generalization of the addition 
formula of theta functions, 
we construct a class of theta functions which are closed 
with respect to the Moyal star product 
of a fixed noncommutative parameter. 
These theta functions can be regarded as bases of the space of 
holomorphic 
homomorphisms between holomorphic line bundles over noncommutative
complex tori. 
\end{abstract}

\maketitle
%
%

\tableofcontents

\section{Introduction}
\noindent
Theta functions are 
associated with various algebraic relations. 
One of them is the addition formula, which also appears 
in the context of the homological mirror symmetry \cite{mirror} 
for elliptic curves \cite{mirror,PZ}, Abelian varieties 
\cite{Fabelian} and also noncommutative real two tori with complex 
structures \cite{foliation,PoSc,KimKim,nchms}. 
It is known that the bases of the space of sections 
for a line bundle on an abelian variety 
is described by theta functions. 
However, in the context of homological mirror symmetry, 
theta functions are regarded rather as the bases of the space of 
holomorphic homomorphisms between two line bundles. 
The composition of two holomorphic homomorphisms is just the 
product of two theta functions, which by the addition formula 
turns out to be a linear combination of theta functions. 
Homological mirror symmetry then asserts that 
such formulas can be reproduced in a geometric way 
by the mirror dual symplectic torus (see subsection \ref{ssec:lagrangian}). 

A noncommutative extension of these stories is given 
in the case of elliptic curves \cite{foliation,PoSc,KimKim,nchms}
based on A.~Schwarz's framework of noncommutative complex tori 
\cite{Stheta,Stensor}. 
However, the conclusion is that the structure constants of the 
product are independent of the noncommutative parameter $\theta$, 
which implies 
that the derived category of holomorphic vector bundles on 
a noncommutative real two-torus is independent of $\theta$ \cite{PoSc}. 

Thus, in order to obtain noncommutative deformations of the 
structure constants, one should discuss 
higher dimensional complex tori. 
In this case, again, an extension of the framework of 
A.~Schwarz's noncommutative complex tori gives various 
explicit noncommutative deformations \cite{ncDG}, which 
includes the deformations described 
in more familiar terminologies 
by the Moyal star product of theta functions. 
In this paper, we present the noncommutative deformation of 
the addition formula of theta functions for 
higher dimensional tori (Theorem \ref{thm:main}). 
For more categorical set-up describing this phenomena, see \cite{ncDG}. 
To explore geometric interpretations 
of this theorem from the mirror dual side should be 
especially interesting. We hope to discuss on it elsewhere. 

In section \ref{sec:comm-theta}, we begin with the commutative case; 
we present explicitly 
the addition formula of theta functions 
corresponding to the holomorphic line bundles on the $n$-dimensional 
complex torus $T^{2n}:=\C^n/(\Z^n\oplus\sqrt{-1}\Z^n)$. 
In section \ref{sec:product}, we explain various aspects of 
the addition formula. 
Though the readers can move ahead to section \ref{sec:nc-theta} directly, 
this section provides us with interesting and pedagogical 
backgrounds on the product of these theta functions, 
together with an introduction to the approach by 
noncommutative complex tori. 
In subsection \ref{ssec:heisenberg}, we explain the relation of 
these theta functions with the {\em theta vectors} introduced by 
A.~Schwarz \cite{Stheta,Stensor} (see also \cite{DKL}) in the 
framework of (non)commutative complex tori. 
In subsection \ref{ssec:twisted}, 
these theta functions or the theta vectors 
are interpreted in terms of 
holomorphic line bundles on complex tori. 
In subsection \ref{ssec:lagrangian}, 
we give an explicit geometric realization of 
the addition formula in the commutative case 
by the mirror dual symplectic torus 
based on the homological mirror symmetry \cite{mirror}. 
In section \ref{sec:nc-theta}, 
we give a noncommutative generalization of this addition formula 
(Theorem \ref{thm:main}). 
Of course, we can replace the product of the addition formula in the 
commutative case by the Moyal star product. 
However, the result is no longer described by any linear combination 
of the theta functions. 
The important point is that we should and in fact can 
find a class of theta functions 
which are closed with respect to the Moyal star product. 
Finally, an example of these noncommutative theta functions in the case of 
complex two-tori is presented in section \ref{sec:exmp}. 

Throughout this paper, any (graded) vector space stands for the one 
over the field $k=\C$.  

\noindent
{\bf Acknowledgments:}\  
I would like to thank 
A,~Kato, T.~Kawai and K.~Saito 
for valuable discussions and useful comments. 
The author is supported by JSPS Research Fellowships for Young Scientists.

 \section{Commutative theta functions}
\label{sec:comm-theta}

The theta function 
$\vartheta: (\R^n/\Z^n\times\R^n/\Z^n)
\times{\frak H}\times\C^n\to\C$ is defined by 
\begin{equation}
 \vartheta[c_1,c_2](\Omega,z):=\sum_{m\in\Z^n} 
 \exp(\pi\sqrt{-1}(m+c_1)^t\Omega (m+c_1) 
 + 2\pi\sqrt{-1}(m+c_1)^t\cdot (z+c_2))\ , 
 \label{cl-theta}
\end{equation}
where $c_1,c_2\in\R^n/\Z^n$ and 
${\frak H}$ is the Siegel upper half plane, that is, 
the space of 
$\C$ valued $n$ by $n$ symmetric matrices whose imaginary parts 
are positive definite. 
Here, for two symmetric matrices $A_a,A_b\in\Mat_n(\Z)$ 
such that $A_{ab}:=A_b-A_a$ is positive definite, 
we define 
\begin{equation}\label{theta-comm}
 \fe^\mu_{ab}(z)=\frac{1}{\sqrt{\det(A_{ab})}}
 \vartheta[0,-A_{ab}^{-1}\mu](\sqrt{-1}A_{ab}^{-1},z) \ ,\qquad 
 \mu\in\Z^n/A_{ab}\Z^n\ , 
\end{equation}
where $\sharp(\Z^n/A_{ab}\Z^n)=\det(A_{ab})$. 
One obtains the following addition formula: 
\begin{thm}\label{thm:comm}
Given three symmetric matrices $A_a,A_b,A_c\in\Mat_n(\Z)$ 
such that $A_{ab}$, $A_{bc}$ 
are positive definite, 
the following product formula holds: 
\begin{equation*}
 \left(\fe_{ab}^\mu\cdot\fe_{bc}^\nu\right)(z)
 =\sum_{\rho\in\Z^n/A_{ac}\Z^n}
 C_{abc,\rho}^{\mu\nu} \fe_{ac}^\rho(z) \ ,
\end{equation*}
where the structure constant $C_{abc,\rho}^{\mu\nu}\in\C$ is given by 
\begin{equation}\label{str-const-comm}
 C_{abc,\rho}^{\mu\nu}
 =\sum_{u\in\Z^n}\delm{A_{ab}}^\mu_{-u+\rho}\delm{A_{bc}}^\nu_u
 \exp\left(-\pi(u-A_{bc}A_{ac}^{-1}\rho)^t
 (A_{ab}^{-1}+A_{bc}^{-1})(u-A_{bc}A_{ac}^{-1}\rho)\right)\ .
\end{equation}
\end{thm}
As explained in the next section, 
in particular, in subsection \ref{ssec:twisted}, 
the collection of these theta functions 
$\{\fe_{ab}^\mu\}_{\mu\in\Z^n/A_{ab}\Z^n}$ 
can be interpreted as the basis of holomorphic 
homomorphisms between a line bundle specified by $A_a$ and 
the one specified by $A_b$ on the $n$-dimensional complex torus 
$T^{2n}=\C^n/(\Z^n+\sqrt{-1}\Z^n)$. 
The addition formula above is then 
interpreted as the composition of the holomorphic homomorphisms. 

Let $\Ob:=\{a,b,\cdots\}$ be a finite collection of labels, 
where any $a\in\Ob$ is associated with a nondegenerate 
symmetric matrix $A_a\in\Mat_n(\Z)$ such that, 
for any $a,b\in\Ob$, $A_{ab}$ is nondegenerate if $a\ne b$. 
For any $a,b\in\Ob$, define a vector space $H^0(a,b)$ over $\C$ as
follows: 
\begin{itemize}
 \item 
If $A_{ab}$ is positive definite, $H^0(a,b)$ is the 
$\det(A_{ab})$-dimensional vector space spanned by the theta functions
$\{\fe_{ab}^\mu\}$. 

 \item If $a=b$, then $H^0(a,b):=\C$. 

 \item If otherwise, then we set $H^0(a,b)=0$. 
\end{itemize}
For any $a,b\in\Ob$, $\Hom(a,a)$ and $\Hom(b,b)$ act on 
$\Hom(a,b)$ from the left and the right, respectively, 
as the trivial multiplication by complex numbers. 
Then, the product formula in Theorem \ref{thm:comm} defines 
an algebraic structure on $\oplus_{a,b\in\Ob}H^0(a,b)$. 
This can in fact be described 
by the zero-th cohomology of an appropriate differential graded category 
(see \cite{ncDG}). 

The main result of this paper 
is a noncommutative extension of Theorem \ref{thm:comm} 
by the Moyal star product (Theorem \ref{thm:main}). 

For the proof of Theorem \ref{thm:comm}, it is convenient to 
prepare the following notion. 
\begin{defn}
Given two symmetric matrices $A_a,A_b\in\Mat_n(\Z)$ such that 
$A_{ab}$ is nondegenerate, 
let $\mu$ be an element in $\Z^n/A_{ab}\Z^n$ and 
$T_{A_{ab}}^\mu:\cS(\R^n)\to C^\infty(T^n)$ 
a linear map defined by 
\begin{equation*}
 (T_{A_{ab}}^\mu\xi)(x)=\sum_{w\in\Z^n}\xi(x+w-A_{ab}^{-1}\mu)\ ,
 \qquad x\in\R^n \ .
\end{equation*}
Here, $\cS(\R^n)$ is the Schwartz space, that is, the space of
functions on $\R^n$ which tend to zero faster than any power of 
$|x|$, $x\in\R^n$. 
 \label{defn:mapT}
\end{defn}
\begin{lem}Let $A_a,A_b,A_c\in\Mat_n(\Z)$ be symmetric matrices such
that $A_{ab}$, $A_{bc}$ and $A_{ac}$ 
are nondegenerate. 
For $\xi_{ab},\xi_{bc}\in\cS(\R^n)$, the following formula holds: 
\begin{equation*}
 (T_{A_{ab}}^\mu \xi_{ab})\cdot(T_{A_{bc}}^\nu \xi_{bc})
 =\sum_{\rho\in\Z^n/A_{ac}\Z^n} (T_{A_{ac}}^\rho \xi_{ac}^\rho)\ ,
\end{equation*}
where $\xi_{ac}^\rho\in \cS(\R^n)$ is defined by 
\begin{equation}\label{tensor-0}
 \xi_{ac}^\rho(x)
 :=\sum_{u\in\Z^n}\delm{A_{ab}}^\mu_{-u+\rho}\delm{A_{bc}}^\nu_u
  \xi_{ab}(x+A_{ab}^{-1}(u-A_{bc}A_{ac}^{-1}\rho))\cdot
  \xi_{bc}(x-A_{bc}^{-1}(u-A_{bc}A_{ac}^{-1}\rho))\ . 
\end{equation}
Here, $\delm{A_{ab}}^\mu_\rho$ is the Kronecker's delta mod 
$\Z^n/A_{ab}\Z^n$, that is, 
\begin{equation*}
 \delm{A_{ab}}^\mu_\rho=
 \begin{cases}
  1 & \qquad \rho-\mu\in A_{ab}\Z^n\ , \\ 
  0 & \qquad \text{otherwise}\ . 
 \end{cases}
\end{equation*}
 \label{lem:Tprod}
\end{lem}
\begin{pf}
By direct calculation, the left hand side is  
\begin{equation*}
 (T_{A_{ab}}^\mu \xi_{ab})\cdot(T_{A_{bc}}^\nu \xi_{bc})(x)
 =\sum_{v\in\Z^n}\delm{A_{ab}}^\mu_{-v} \xi_{ab}(x+A_{ab}^{-1}v)
  \sum_{v'\in\Z^n}\delm{A_{bc}}^\nu_{-v'} \xi_{bc}(x+A_{bc}^{-1}v')\ .
\end{equation*}
By the transformation 
\begin{equation*}
 \bp v \\ v'\ep = 
 \bp \1_n & A_{ab} \\ -\1_n & A_{bc} \ep 
 \bp u \\ w \ep - \bp \rho \\ \0_n \ep\ ,
\end{equation*}
the equation above is rewritten as 
\begin{equation*}
  (T_{A_{ab}}^\mu \xi_{ab})\cdot(T_{A_{bc}}^\nu\xi_{bc})(x)
 =\sum_{\rho\in\Z^n/A_{ac}\Z^n}
 \sum_{u,w\in\Z^n}\delm{A_{ab}}^\mu_{-u+\rho}\delm{A_{bc}}^\nu_u 
 \xi_{ab}(x+w+A_{ab}^{-1}(u-\rho))\xi_{bc}(x+w-A_{bc}^{-1}u)\ .
\end{equation*}
On the other hand, the right hand side can be computed directly as 
\begin{equation*}
 \begin{split}
 (T_{A_{ac}}^\rho \xi_{ac}^\rho)(x)
 &=\sum_{\rho\in\Z^n/A_{ac}\Z^n}\sum_{u,w\in\Z^n}
 \delm{A_{ab}}^\mu_{-u+\rho}\delm{A_{bc}}^\nu_u \\
 &\qquad 
 \xi_{ab}(x+w-A_{ac}^{-1}\rho+A_{ab}^{-1}(u-A_{bc}A_{ac}^{-1}\rho))\cdot
  \xi_{bc}(x+w-A_{ac}^{-1}\rho-A_{bc}^{-1}(u-A_{bc}A_{ac}^{-1}\rho))\\
 &=\sum_{u\in\Z^n}\delm{A_{ab}}^\mu_{-u+\rho}\delm{A_{bc}}^\nu_u
  \xi_{ab}(x+w+A_{ab}^{-1}(u-\rho))\cdot
  \xi_{bc}(x+w-A_{bc}^{-1}u)\ .
 \end{split}
\end{equation*}
Thus, the left hand side coincides with the right hand side. 
\qed\end{pf}
For $A_a,A_b$ such that $A_{ab}\in\Mat_n(\Z)$ is positive definite, 
define a function $e_{ab}\in\cS(\R^n)$ by 
\begin{equation}\label{theta-vector-comm1}
 e_{ab}(x)=\exp\left(-\pi x^tA_{ab}x\right)\ .
\end{equation}
Then, by the Poisson resummation formula (see \cite{mumford1}, p195-197), 
one can rewrite the theta functions $\fe_{ab}^\mu(z)$ as 
\begin{equation}\label{e_ab-comm}
 \fe_{ab}^\mu(z):=T_{A_{ab}}^\mu(e_{ab})(z)\ ,
 \qquad \mu\in\Z^n/A_{ab}\Z^n\ ,
\end{equation}
where, for $T_{A_{ab}}^\mu(e_{ab})\in\cS(\R^n)$, 
$T_{A_{ab}}^\mu(e_{ab})(z)$ stands for the holomorphic extension. 

Thus, for symmetric matrices $A_a,A_b,A_c\in\Mat_n(\Z)$ 
such that $A_{ab}$ and $A_{bc}$ are positive definite, 
apply Lemma \ref{lem:Tprod} with $\xi_{ab}=e_{ab}$, $\xi_{bc}=e_{bc}$, 
and the holomorphic extension leads to Theorem \ref{thm:comm}.

 \section{Various interpretations of the product formula}
\label{sec:product}

In this section, we give various interpretations of Theorem \ref{thm:comm}. 

 \subsection{The tensor product of Heisenberg modules}
\label{ssec:heisenberg}

Theorem \ref{thm:comm} can be understood directly 
in A.~Schwarz's 
framework of noncommutative complex tori \cite{Stheta,Stensor}. 
A {\em noncommutative torus} $\cA^{d}_\theta$ is an algebra 
defined by unitary generators 
$U_1,\cdots,U_d$ with relations 
\begin{equation}\label{nctorus}
 U_iU_j=e^{-2\pi\sqrt{-1}\theta_{ij}}U_jU_i\ ,\qquad
 \theta_{ij}=-\theta_{ji}\in\R\ 
\end{equation}
for $i,j=1,\cdots,d$. 
Now, we shall consider $2n$-dimensional commutative torus 
$\cA^{2n}:=\cA^{2n}_{\theta=0}$. 
Namely, $\cA^{2n}$ is thought of as the space of functions 
on a $2n$-dimensional commutative torus $T^{2n}$. 
Thus, the generators $U_1,\cdots,U_{2n}$ now commute with each other. 

A pair $E_a:=(E_{A_a},\nabla_a)$ 
of a finitely generated projective module $E_{A_a}$, 
called a {\em Heisenberg module} (see \cite{KS}), 
with a constant curvature connection $\nabla_a$ 
is constructed as follows. 
The Heisenberg module is defined by 
$$E_{A_a}:=\cS(\R^n\times(\Z^n/A_{a}\Z^n))$$ 
for a fixed nondegenerate symmetric matrix $A_a\in\Mat_n(\Z)$. 
The right action of $\cA^{2n}$ on $E_{A_a}$ is defined by 
specifying the right action of each generator; 
for $\xi_a\in E_{A_a}$, it is given by 
\begin{equation}\label{Uaction}
 \begin{split}
 (U_i\xi_a)(x;\mu)& 
 = e^{2\pi\sqrt{-1}(x_i+(A_{a}^{-1}\mu)_i))} \xi_a(x;\mu)\ ,\\
 (U_{n+i}\xi_a)(x;\mu)& =\xi_a(x+A_{a}^{-1}t_i;\mu-t_i)
 \ ,\qquad i=1,\cdots, n \ ,
 \end{split}
\end{equation}
where $x:=(x_1\cdots x_n)^t\in\R^n$ ($^t$ indicates the transpose), 
$\mu\in\Z^n/A_{a}\Z^n$ 
and $t_i\in\R^n$ is defined by $(t_1 \cdots t_n)=\1_n$. 
A constant curvature connection 
$\nabla_{a,i}:E_{A_a}\to E_{A_a}$, $i=1,\cdots,2n$, 
is given by 
\begin{equation}\label{ccc}
 (\nabla_{a,1}\cdots \nabla_{a,2n})^t= 
 \bp
  \1_n &  \\ 
    & -A_{a} 
 \ep
 \bp \partial_x \\ 2\pi\sqrt{-1}x \ep\ ,
\end{equation}
where 
$\partial_x:=\left(\bps\fpartial{x_1}& \cdots &
\fpartial{x_n}\eps\right)^t$, 
whose curvature 
$F_a:=\{\frac{\sqrt{-1}}{2\pi}[\nabla_{a,i},\nabla_{a,j}]\}_{i,j=1,\cdots,2n}$ 
is 
\begin{equation*}
 F_a:=\bp \0_n & A_a \\ -A_a & \0_n\ep\ .
\end{equation*}
The generators of the endomorphism algebra is the same as 
$U_i$, $i=1,\cdots,2n$: 
\begin{equation*}
 \begin{split}
 (\xi_aZ_i)(x;\mu)& 
 =\xi_a(x;\mu)\, e^{2\pi\sqrt{-1}(x_i+(A_{a}^{-1}\mu)_i))}\ ,\\
 (\xi_aZ_{n+i})(x;\mu)& =\xi_a(x+A_{a}^{-1}t_i;\mu-t_i)
 \ ,\qquad i=1,\cdots, n \ .
 \end{split}
\end{equation*}
Namely, the endomorphism algebra also forms 
a commutative torus $\cA^{2n}$. 

Given $E_a$ and $E_b$ such that $A_{ab}$ is nondegenerate, 
the space $\Hom(E_a,E_b)$ is defined 
again as the Schwartz space
$\Hom(E_a,E_b):=\cS(\R^n\times(\Z^n/A_{ab}\Z^n))$. 
For $\xi_{ab}\in\Hom(E_a,E_b)$, the right action of $\cA^{2n}$, 
generated by $U_i$, $i=1,\cdots,2n$, 
and the left action of $\cA^{2n}$, 
generated by $Z_i$, $i=1,\cdots,2n$, are defined by 
\begin{equation*}
  \begin{split}
 & (U_i\xi_{ab})(x;\mu) = 
 e^{2\pi\sqrt{-1}(x_i+(A_{ab}^{-1}\mu)_i))}\, \xi_{ab}(x;\mu)\ ,\qquad
 (U_{n+i}\xi_{ab})(x;\mu) =\xi_{ab}(x+A_{ab}^{-1}t_i;\mu-t_i)\ ,\\
 &(\xi_{ab}Z_i)(x;\mu) =\xi_{ab}(x;\mu)\, 
 e^{2\pi\sqrt{-1}(x_i+(A_{ab}^{-1}\mu)_i))} \ ,\qquad 
 (\xi_{ab}Z_{n+i})(x;\mu) =\xi_{ab}(x+A_{ab}^{-1}t_i;\mu-t_i)\ ,
 \end{split}
\end{equation*}
where $\mu\in\Z^n/A_{ab}\Z^n$. 
In fact, all these generators $U_i$ and $Z_i$, $i=1,\cdots,2n$, 
commute with each other. 
The constant curvature connection 
$\nabla_i:\Hom(E_a,E_b)\to\Hom(E_a,E_b)$, 
$i=1,\cdots, 2n$, is given by 
\begin{equation*}
 (\nabla_1 \cdots \nabla_{2n})^t := 
 \bp
  \1_n &  \\ 
    & -A_{ab} 
 \ep
 \bp \partial_x \\ 2\pi\sqrt{-1}x \ep\ .
\end{equation*}
For $\xi_{ab}\in\Hom(E_a,E_b)$ and 
$\xi_{bc}\in\Hom(E_b,E_c)$, 
the tensor product 
$m:\Hom(E_a,E_b)\otimes\Hom(E_b,E_c)\to\Hom(E_a,E_c)$ 
is defined by 
\begin{equation} \label{high-tensor}
 m(\xi_{ab},\xi_{bc})(x,\rho)
 =\sum_{u\in\Z^n}
  \xi_{ab}(x+A_{ab}^{-1}(u-A_{bc}A_{ac}^{-1}\rho),-u+\rho)\cdot
  \xi_{bc}(x-A_{bc}^{-1}(u-A_{bc}A_{ac}^{-1}\rho),u)\ .
\end{equation}
One can see that this tensor product formula is 
just the definition of $\xi_{ac}^\rho$ 
in eq.(\ref{tensor-0}). 
This tensor product is in fact associative and 
the connection $\nabla_i:\Hom(E_a,E_b)\to\Hom(E_a,E_b)$ 
satisfies the Leibniz rule with respect to this product 
(see \cite{ncDG}).
\footnote{In \cite{ncDG}, left modules in this paper is flipped 
to be right modules. 
The relation of the conventions between this paper and \cite{ncDG} 
is as follows. 
First, consider a bimodule $\Hom(E_a,E_b)$ in this paper. 
Replace $A_a$ by $-A_b$ and $A_b$ by $-A_a$. 
Then, one gets a bimodule in \cite{ncDG}.  
In both cases, a left/right module $E_{A_b}$ is obtained by 
setting $A_a=0$. }

Now suppose we consider a $n$-dimensional {\em complex} torus 
$T^{2n}:=\C^n/(\Z^n\oplus\sqrt{-1}\Z^n)$. 
For $E_a=(E_{A_a},\nabla_a)$ a Heisenberg module 
with the constant curvature connection, 
the holomorphic structure 
$\nabb_{a,i}:E_{A_a}\to E_{A_a}$, $i=1,\cdots, n$, is defined by 
\begin{equation*}
 \nabb_{a,i}=\nabla_{a,i}+\sqrt{-1}\nabla_{a,n+i}\ .
\end{equation*}
Also, for given $E_a, E_b$, the holomorphic structure 
$\nabb_i:\Hom(E_a,E_b)\to\Hom(E_a,E_b)$, 
$i=1,\cdots,n$, is defined in the same way: 
\begin{equation*}
 \nabb_i:=\nabla_i+\sqrt{-1}\nabla_{n+i}\ ,\qquad i=1,\cdots, n\ .
\end{equation*}
When $A_{ab}$ is positive definite, 
the space $H^0(E_a,E_b):=
\cap_{i=1}^n\Ker(\nabb_i:\Hom(E_a,E_b)\to\Hom(E_a,E_b))$ 
forms a $\det(A_{ab})$-dimensional vector space. 
The bases $e_{ab}^\mu$, $\mu\in\Z^n/A_{ab}\Z^n$, 
are called A.~Schwarz's {\em theta vectors} \cite{Stheta} 
(see also \cite{DKL}), which are just the function
$e_{ab}\in\cS(\R^n)$ defined in eq.(\ref{theta-vector-comm1}): 
\begin{equation}\label{theta-vector}
 e_{ab}^\mu(x,\rho)=\delm{A_{ab}}^\mu_\rho
 \exp\left(-\pi x^tA_{ab}x\right)\ . 
\end{equation}
The Leibniz rule of $\nabb$ 
then guarantees that the tensor product $m(e_{ab}^\mu,e_{bc}^\nu)$ 
turns out to be the linear combination of $e_{ac}^\rho$, 
$\rho\in\Z^n/A_{ac}\Z^n$. 

This approach by Heisenberg modules allows us various 
noncommutative deformations of these structures (see \cite{ncDG}), 
but some of such deformations 
can be lifted to theta functions 
as the Moyal star product; 
the consequence is the one presented in section \ref{sec:nc-theta}.

\subsection{Holomorphic line bundles on tori}
\label{ssec:twisted}

In this subsection, 
the theta functions $\{\fe_{ab}^\mu\}$ 
in eq.(\ref{e_ab-comm}), 
or equivalently, the theta vectors $\{e_{ab}^\mu\}$ 
in eq.(\ref{theta-vector}), 
are interpreted in terms of 
holomorphic line bundles on complex tori.

Given a $d$-dimensional-torus $T^{d}=\R^{d}/\Lambda$, 
$\Lambda:=\Z^{d}$, 
let $\pi: \R^{d}\to \R^{d}/\Z^{d}$ be the projection. 
A vector bundle $p:E\to T^{d}$ is constructed as the 
pullback $\pi^*E$ together with the action of $\Lambda$, 
where 
\begin{equation*}
 \pi^*E=\{ (x,\xi)\in \R^{2n}\times E\ |\ \pi(x)= p(\xi) \}\ .
\end{equation*}
When $p:E\to T^d$ is a rank $q$ vector bundle, $\pi^*E$ is a 
rank $q$ trivial vector bundle over $\R^{d}$. 
The $\Lambda\ni\lambda$ action on the sections of 
$\pi^*E\simeq\R^{d}\times\C^q$ is defined by 
\begin{equation}
 \xi(x+\lambda):=c_\lambda(x)\xi(x)\ ,\qquad 
 \xi\in \Gamma(\pi^*E)\simeq(C^\infty(\R^d))^{\oplus q}\ ,\qquad 
 c_\lambda\in U(q;C^\infty(\R^d))\ .
 \label{xi-twist}
\end{equation}
We require that this $c$ satisfies the following condition: 
\begin{equation}
 c_{\lambda'}(x+\lambda)c_\lambda(x)=
 c_{\lambda+\lambda'}(x)\ .
 \label{cocycle-cd}
\end{equation}
Thus, $c_\gamma$ is regarded as a transition function 
of the vector bundle $E$. 
A connection $\nabla_i:\Gamma(\pi^*E)\to\Gamma(\pi^*E)$, 
$i=1,\cdots, d$, is defined so that 
the following compatibility conditions hold: 
\begin{equation}
 (\nabla_i)(x+\lambda)
=c_\lambda(x)(\nabla_i)(x)c^{-1}_\lambda(x)\ ,
 \label{conn-compatible}
\end{equation} 
where the {\em curvature} is defined by 
\begin{equation*}
 F=\{F_{ij}\}_{i,j=1,\cdots,d}\ ,\qquad 
 F_{ij}:=\frac{\sqrt{-1}}{2\pi}[\nabla_i,\nabla_j]\ .
\end{equation*}
Now, let us consider a complex torus
$T^{2n}:=\C^n/(\Z^n\oplus\sqrt{-1}\Z^n)$, 
where we denote the coordinates of the covering space $\C^n$ 
by $z:=(z_1\cdots z_n)^t$, 
$z_i:=x_i+\sqrt{-1} y_i$, $i=1,\cdots,n$. 
For a nondegenerate symmetric matrix $A_a\in\Mat_n(\Z)$, 
the space of sections $\ti{E}_{A_a}$ 
of a line bundle ($q=1$ case) on $T^{2n}$ 
is constructed by setting 
\begin{equation*}
 c_{(\lambda_x,0)}(x,y)=1\ ,\qquad 
 c_{(0,\lambda_y)}(x,y)=e^{-2\pi\sqrt{-1}x^t A_a\lambda_y}\cdot 1\ ,
\end{equation*}
where $x:=(x_1\cdots x_n)^t$, $y:=(y_1\cdots y_n)^t$ and 
$\lambda_x,\lambda_y\in\Z^n$ such that 
$\lambda=(\lambda_x,\lambda_y)\in\Lambda$. 
In order to show that this transition function $c_\gamma$
satisfies the condition (\ref{cocycle-cd}), 
it is enough to check 
\begin{equation*}
 c_{(\lambda_x,0)}^{-1}(x+\lambda_x,y)
 c^{-1}_{(0,\lambda_y)}(x,y+\lambda_y)
 c_{(\lambda_x,0)}(x,y+\lambda_y)c_{(0,\lambda_y)}(x,y)=1\ .  
\end{equation*}
The general form of sections in $\ti{E}_{A_a}$ is 
described as a function on the covering space $\R^{2n}$ 
with coordinates $(x,y)$ satisfying eq.(\ref{xi-twist}); 
it is given as a natural extension of 
the two dimensional case (\cite{GRT, Ho, MZ} and see \cite{KS}, 
the vector bundles constructed there are called {\em twisted bundles}): 
\begin{equation*}
\ti{\xi}_a(x,y)=\sum_{w\in\Z^n}\sum_{\mu\in\Z^n/A\Z^n}
\!\exp\left(2\pi\sqrt{-1} y^t\left(
-A_a\left(x+w\right)+\mu\right)\right)
\xi^\mu_a\left(x+w-A_a^{-1}\mu\right)\, ,\ \ \xi^\mu_a\in\cS(\R^n)\, .
\end{equation*}
For $\xi^\mu_a(x)=:\xi_a(x,\mu)$, 
$\xi_a\in\cS(\R^n\otimes(\Z^n/A_a\Z^n))=E_{A_a}$, 
we regard $\ \ti{}\ $ in the formula above as the isomorphism 
from $E_{A_a}$ to $\ti{E}_{A_a}$ 
which sends $\xi_a$ to $\ti{\xi}_a$. 
This line bundle can be equipped with the following constant
curvature connection 
$\{\nabla_{a,i}:\ti{E}_{A_a}\to\ti{E}_{A_a}\}_{i=1,\cdots,2n}$ 
with its curvature $F_a$: 
\begin{equation*}
 (\nabla_{a,1},\cdots,\nabla_{a,n})^t
 =\partial_x + 2\pi\sqrt{-1}A y\ ,\quad  
 (\nabla_{a,n+1},\cdots,\nabla_{a,2n})^t=\partial_y\ ,\qquad 
 F_a= \bp \0_n & A_a \\ -A_a & \0_n \ep\ ,
 \label{gf}
\end{equation*}
where $\partial_{a,x}:=(\fpartial{x_1}\cdots\fpartial{x_n})^t$, 
$\partial_{a,y}:=(\fpartial{y_1}\cdots\fpartial{y_n})^t$.  
Let us define the generators of the space $C^\infty(T^{2n})$ 
of functions by 
\begin{equation*}
 \ti{U}_i=e^{\pi\sqrt{-1} x_i}\ ,\quad 
 \ti{U}_{n+i}=e^{\pi\sqrt{-1} y_i}\ ,\qquad i=1,\cdots, n\ . 
\end{equation*}
Then, the relationship of $\ti{E}_a:=(\ti{E}_{A_a},\nabla_a)$ 
with $E_a=(E_{A_a},\nabla_a)$ in the previous subsection 
can be summarized as follows: 
for $\xi_{a}\in E_{A_a}$, 
\begin{equation*}
 \ti{U}_i\ti{\xi}_a=\wti{U_i\xi}_a\ ,\quad 
 \ti{U}_{n+i}\ti{\xi}_a=\wti{U_{n+i}\xi_a}\ ,\quad 
 \nabla_{a,i}\ti{\xi}_a=\wti{\nabla_{a,i}\xi}\ ,\quad 
 \nabla_{a,n+i}\ti{\xi}_a=\wti{\nabla_{a,n+i}\xi_a}\ , 
 \qquad i=1,\cdots, n\ .
\end{equation*}
In a similar way, for given $\ti{E}_a$ and $\ti{E}_b$ 
such that $A_{ab}$ is nondegenerate, 
the space $\Hom(\ti{E}_a,\ti{E}_b)$ of homomorphisms 
from $\ti{E}_a$ from $\ti{E}_b$ is the space 
whose elements are described of the form: 
\begin{equation}\label{general}
 \ti{\xi}_{ab}(x,y)=\sum_{w\in\Z^n}\sum_{\mu\in\Z^n/A_{ab}\Z^n}
 \!\exp\left(2\pi\sqrt{-1} y^t\left(
 -A_{ab}\left(x+w\right)+\mu\right)\right)
 \xi_{ab}^\mu\left(x+w-A_{ab}^{-1}\mu\right)\, ,
\end{equation}
for $\xi_{ab}^\mu\in\cS(\R^n)$, 
where the compatible constant curvature connection 
$\nabla_i:\Hom(\ti{E}_a,\ti{E}_b)\to\Hom(\ti{E}_a,\ti{E}_b)$, 
$i=1,\cdots, n$, is given by 
\begin{equation*}
 (\nabla_{1},\cdots,\nabla_{n})^t
 :=\partial_x + 2\pi\sqrt{-1}A_{ab} y\ ,\quad 
(\nabla_{n+1},\cdots,\nabla_{2n})^t 
 :=\partial_y\ ,\qquad 
 F_{ab}= \bp \0_n & A_{ab} \\ -A_{ab} & \0_n \ep\ .
\end{equation*}
Again, 
for $\xi_{ab}^\mu(x)=:\xi_{ab}(x,\mu)$, 
$\xi_{ab}\in\cS(\R^n\otimes(\Z^n/A_{ab}\Z^n))=\Hom(E_a,E_b)$, 
$\ \ti{}\ $ in eq.(\ref{general}) is regarded as the isomorphism 
from $\Hom(E_a,E_b)$ to $\Hom(\ti{E}_a,\ti{E}_b)$ 
which sends $\xi_{ab}$ to $\ti{\xi}_{ab}$. 

Actually, for $E_a,E_b,E_c$, 
$\xi_{ab}\in\Hom(E_a,E_b), \xi_{bc}\in\Hom(E_b,E_c)$ 
and the corresponding elements 
$\ti{\xi}_{ab}\in\Hom(\ti{E}_a,\ti{E}_b),
\ti{\xi}_{bc}\in\Hom(\ti{E}_b,\ti{E}_c)$, 
the pointwise product $\ti{\xi}_{ab}\cdot\ti{\xi}_{bc}$ turns out to be 
\begin{equation*}
\ti{\xi}_{ab}\cdot\ti{\xi}_{bc}=\wti{m(\xi_{ab},\xi_{bc})}\ ,
\end{equation*}
where $m$ is the tensor product of the Heisenberg modules defined 
in eq.(\ref{high-tensor}). 
The proof is essentially the same as that of Lemma
\ref{lem:Tprod}. 

Now, for $T^{2n}$ as a {\em complex} torus, 
the holomorphic structure 
$\{\nabb_{a,i}:\ti{E}_a\to\ti{E}_a\}_{i=1,\cdots,n}$ 
is defined by 
$\nabb_{a,i}:=\nabla_{a,i}+\sqrt{-1}\nabla_{a,n+i}$. 
Similarly, given $\ti{E}_a$ and $\ti{E}_b$, the holomorphic structure 
$\{\nabb_{i}:\Hom(\ti{E}_a,\ti{E}_b)
\to\Hom(\ti{E}_a,\ti{E}_b)\}_{n=1,\cdots,n}$ 
is defined by $\nabb_{i}:=\nabla_{i}+\sqrt{-1}\nabla_{n+i}$. 
The space of holomorphic sections in $\Hom(\ti{E}_a,\ti{E}_b)$ 
is then defined by 
$H^0(\ti{E}_a,\ti{E}_b):=
\cap_{i=1}^n
\Ker(\nabb_i:\Hom(\ti{E}_a,\ti{E}_b)\to\Hom(\ti{E}_a,\ti{E}_b))$. 
This space $H^0(\ti{E}_a,\ti{E}_b)$ 
forms a $\det(A_{ab})$-dimensional vector
space spanned by $\{\ti{e}_{ab}^\mu\}$, the extension of the theta vectors 
$\{e_{ab}^\mu\}_{\mu\in\Z^n/A_{ab}\Z^n}$ in (\ref{theta-vector})  
by eq.(\ref{general}). 
Also, the explicit relation of these $\ti{e}_{ab}^\mu$ with 
the theta functions $\fe_{ab}^\mu$ (\ref{e_ab-comm}) is given by 
\begin{equation*}
 \fe_{ab}^\mu(z)=\exp\left(\pi y^tA_{ab}y\right)\cdot\ti{e}_{ab}^\mu(x,y)\ .
\end{equation*}

 \subsection{Lagrangian submanifolds and triangles}
\label{ssec:lagrangian}

The homological mirror symmetry \cite{mirror} asserts 
that the product $m(\fe_{ab}^\mu,\fe_{bc}^\nu)$ 
can also be derived from geometry of 
the mirror dual torus $\Th^{2n}$, 
a symplectic $2n$-dimensional torus with the symplectic structure 
\begin{equation}\label{omega}
 \omega=\bp \0_n & -\1_n \\ \1_n & \0_n \ep \ .
\end{equation}
For the covering space $\R^{2n}$ of $\Th^{2n}$, 
let $\pi:\R^{2n}\to\Th^{2n}$ be the natural projection. 
The coordinates for $\R^{2n}$ is denoted 
$(x_1,\cdots,x_n,\yh_1,\cdots,\yh_n)$

The affine lagrangian submanifold mirror dual to 
$E_a=(E_{A_a},\nabla_a)$ over $\cA^{2n}$, the space of functions on 
$T^{2n}$, 
is defined by the image of the affine subspace in $\R^{2n}$ 
\begin{equation*}
 L_a: \yh=A_a x
\end{equation*}
by the projection $\pi:\R^{2n}\to\Th^{2n}$. 
Thus, we have 
\begin{equation*}
\pi^{-1}\pi(L_a)=\{\yh=A_a x+c_a,\ c_a\in\Z\}\ .
\end{equation*}
Let us define the space of morphisms 
$\Hom(L_a,L_b)$ which is isomorphic to the 
$\sharp(\Z^n/A_{ab}\Z^n)$-dimensional vector space 
$H^0(E_a,E_b)$ in subsection \ref{ssec:heisenberg}. 
Denote the basis of $\Hom(L_a,L_b)$ by 
$v_{ab}^\mu$, $\mu\in\Z^n/A_{ab}\Z^n$, 
to which is associated the image of the 
intersection point of $\yh=A_b x+\mu$ with $\yh=A_a x$ in $\C^n$ 
by $\pi:\C^n\to\Th^{2n}$. 
One can see that actually the intersection point of 
$\pi(L_a)$ and $\pi(L_b)$ in $\Th^{2n}$ is 
$\sharp(\Z^n/A_{ab}\Z^n)=\det(A_{ab})$. 
For a bases $v_{ab}$ of $\Hom(L_a,L_b)$, we denote the corresponding 
point in $\Th^{2n}$ also by $v_{ab}$, which defines the set 
$\ti{V}_{ab}:=\pi^{-1}(v_{ab})$ of points 
in the covering space $\R^{2n}$. 

The structure constant $C_{abc,\rho}^{\mu\nu}\in\C$ 
(\ref{str-const-comm}) can be identified 
with the sum of the exponentials of the symplectic areas 
of the triangles 
$\ti{v}_{ab}\ti{v}_{bc}\ti{v}_{ac}$ for any 
$\ti{v}_{ab}\in\ti{V}_{ab}$, $\ti{v}_{bc}\in\ti{V}_{bc}$ and 
$\ti{v}_{ac}\in\ti{V}_{ac}$ 
with respect to the symplectic structure $\omega$ in
eq.(\ref{omega}), where the triangles related by parallel 
transformations on the covering space $\R^{2n}$ are 
identified with each other. 

It is calculated as follows. 
Consider three affine subspaces $L_a',L_b',L_c'$ in $\R^2$ as follows: 
\begin{equation*}
 L_a': \yh = A_a x + c_a\ ,\quad 
 L_b': \yh = A_b x + c_b\ ,\quad 
 L_c': \yh = A_c x + c_c\ . 
\end{equation*}
If $A_{ab}$ is nondegenerate, 
the intersection of $L_a'$ and $L_b'$ is a point $v_{ab}$; the 
coordinates $\left(\bps x\\ \yh\eps\right)$ are: 
\begin{equation*}
 v_{ab}=
 \bp
 -(A_{ab})^{-1}(c_b-c_a)\\
 -A_aA_{ab}^{-1}c_b +A_bA_{ab}^{-1}c_a
 \ep
\ .
\end{equation*}
Now, assume that $A_{ab}$ and $A_{bc}$ are 
positive definite. Then, $A_{ac}$ is also positive definite. 
The three intersection points $v_{ab}, v_{bc}, v_{ac}$ 
form a triangle, where the edges $(v_{ab}v_{bc})$, $(v_{bc}v_{ac})$, 
$(v_{ac}v_{ab})$ belong to $L_b'$, $L_c'$ and $L_a'$, respectively. 
The symplectic area of the triangle is defined by 
\begin{equation*}
 (v_{ab}-v_{ac})^t \omega (v_{bc}-v_{ac})
 =\bp (c_c-c_a)^t & (c_b-c_a)^t \ep 
 \bp A_{bc}^{-1} & A_{ac}^{-1} \\
    A_{ab}^{-1}A_{ac}A_{bc}^{-1} & A_{ab}^{-1}\ep
 \bp c_b-c_c \\ c_a-c_c\ep \ .
\end{equation*}
Let us put $c_a=0$, $c_b=u'$ and $c_c=-\rho$ so that 
$\pi(v_{ac})=v_{ac}^\rho$. 
Then, consider 
\begin{equation*}
 \sum_{u'}\delm{A_{ab}}^\mu_{-u'}\delm{A_{bc}}^\nu_{u'+\rho}
 \exp\left( \left(\bps -\rho^t & {u'}^t \eps\right) 
 \left(\bps A_{bc}^{-1} & A_{ac}^{-1} \\
    A_{ab}^{-1}A_{ac}A_{bc}^{-1} & A_{ab}^{-1}\eps\right)
 \left(\bps u'+\rho \\ \rho\eps\right) \right)\ ,
\end{equation*}
where $\delm{A_{ab}}^\mu_{-u'}$ and $\delm{A_{bc}}^\nu_{u'+\rho}$ 
correspond to the condition of $\pi(v_{ab})=v_{ab}^\mu$ and 
$\pi(v_{bc})=v_{bc}^\nu$, respectively. 
One can see that, 
by the replacement $u'+\rho=:u$, 
this coincides with the structure constant $C_{abc,\rho}^{\mu\nu}$ 
of the product of the theta 
functions in eq.(\ref{str-const-comm}).

 \section{Noncommutative theta functions}
\label{sec:nc-theta}

The Moyal star product (\cite{Moyal}) is an associative 
noncommutative product on functions on a flat space. 
It gives the first example of deformation quantization \cite{BFFLS} 
and is also be used as a building block of deformation quantization on 
arbitrary symplectic manifolds (see \cite{OMY91,GR}). 
A Moyal star product on functions on $\C^n$ is defined by
\begin{equation*}
 (f*g)(z)=f(z)
 e^{-\frac{\sqrt{-1}}{4\pi}\lpartial{z}\theta\rpartial{z}}
 g(z)\ , 
\end{equation*}
where 
$\lpartial{z}\theta\rpartial{z}:=\sum_{i,j=1}^n
\flpartial{z^i}\theta^{ij}\frpartial{z^j}$. 
Note that this skewsymmetric matrix $\theta\in\Mat_{n}(\R)$ 
can be thought of as the restriction of the 
$\theta=\{\theta_{ij}\}_{i,j=1,\cdots,2n}$ in eq.(\ref{nctorus}) 
to $\theta=\{\theta_{ij}\}_{i,j=1,\cdots,n}$. 
\footnote{This skewsymmetric matrix $\theta\in\Mat_n(\R)$ corresponds 
to $\theta_1$ in \cite{ncDG}.} 

Now, for two symmetric matrices $A_a,A_b\in\Mat_n(\C)$ 
such that $A_{ab}$ is nondegenerate, 
the following matrix $M_{ab}\in\Mat_n(\C)$, 
\begin{equation*}
  M_{ab}:=
 \left(\1_n+\frac{\sqrt{-1}}{2}A_{ab}^+\theta\right)^{-1} 
 A_{ab}\ , \qquad A_{ab}^+:=A_a+A_b\ ,
\end{equation*}
is symmetric if and only if the 
the following condition holds: 
\begin{equation}\label{DGcd}
 A_a\theta A_a=A_b\theta A_b\ . 
\end{equation}
For $A_a,A_b\in\Mat_n(\C)$ satisfying the condition (\ref{DGcd}), 
the real part of $M_{ab}$ is positive definite 
if and only if $A_{ab}$ is positive definite (see \cite{Ig}, p.5). 
For two matrices $A_a,A_b\in\Mat_n(\C)$ such that 
$A_{ab}$ is positive definite, 
define theta functions $\fe_{ab}^\mu$, $\mu\in\Z^n/A_{ab}\Z^n$, by 
\begin{equation}\label{nc-theta}
 \fe_{ab}^\mu(z)=
 \frac{\det(\1_n+\sqrt{-1} A_a\theta)^{\ov{4}}
\det(\1_n+\sqrt{-1} A_b\theta)^{\ov{4}}}
{\det(A_{ab})^{\ov{2}}}\, 
\vartheta[0,-A_{ab}\mu](\sqrt{-1} M_{ab}^{-1},z)\ .
\end{equation}
It is clear that these theta functions actually coincides with 
those in eq.(\ref{theta-comm}) if $\theta=0$. 

Then, we get the $*$ product formula of these noncommutative theta
functions. 
\begin{thm}\label{thm:main}
For a fixed $\theta$, consider a set of 
symmetric matrices $A_a,A_b,A_c\in\Mat_n(\Z)$ such that 
$A_a\theta A_a=A_b\theta A_b=A_c\theta A_c$ and 
$A_{ab},A_{bc}\in\Mat_n(\C)$ are positive definite. 
Then, the following product formula holds:
\begin{equation*}
 \left(\fe_{ab}^\mu * \fe_{bc}^\nu\right)(z)= \sum_{\rho\in\Z^n/A_{ac}\Z^n} 
 C_{abc,\rho}^{\mu\nu} \fe_{ac}^\rho(z)\ , 
\end{equation*}
\begin{equation*}
 C_{abc,\rho}^{\mu\nu}
 :=\sum_{u\in\Z^n}\delm{A_{ab}}^\mu_{-u+\rho}\delm{A_{bc}}^\nu_u
 \exp\left(-\pi(u-A_{bc}A_{ac}^{-1}\rho)^t
 \left((A_{ab}^{-1}+A_{bc}^{-1})(\1+\sqrt{-1}A_b\theta)^{-1}\right)
 (u-A_{bc}A_{ac}^{-1}\rho)\right)\ .
\end{equation*}
\end{thm}
Note that the matrix 
$(A_{ab}^{-1}+A_{bc}^{-1})(\1+\sqrt{-1}A_b\theta)^{-1}
\in\Mat_n(\C)$ is already symmetric.  
\begin{pf}
Again, by the Poisson resummation formula, 
the theta functions $\{\fe_{ab}^\mu\}$ in eq.(\ref{nc-theta}) 
can be rewritten as $\fe_{ab}^\mu(z)=T_{A_{ab}}^\mu(e_{ab})(z)$, 
where 
\begin{equation*}
  e_{ab}(x):=C_{ab}\cdot e^{-\pi x^t M_{ab} x}\ ,\qquad 
 C_{ab}:=
 \frac{\det(\1_n+\sqrt{-1} A_a\theta)^{\ov{4}}
\det(\1_n+\sqrt{-1} A_b\theta)^{\ov{4}}}
{\det(\1_n+\2i A_{ab}^+\theta)^{\ov{2}}}\ .
\end{equation*}
As in the commutative case in subsection \ref{ssec:heisenberg}, 
one can consider the corresponding Heisenberg modules with 
a constant curvature connection $\nabla$, 
where the tensor product is given just by replacing the product
$\cdot$ in the right hand side of eq.(\ref{tensor-0}) by 
the star product, the constant curvature connection $\nabla$ 
satisfies the Leibniz rule 
with respect to the tensor product, 
and the the theta vectors are obtained just as the function $e_{ab}$ 
above \cite{ncDG}. 
The Leibniz rule of $\nabla$ then guarantees that the tensor 
product $m(e_{ab}^\mu,e_{bc}^\nu)$ is a linear combination of 
$e_{ac}^\rho$. 
The appropriate coefficients $C_{ab}\in\C$ and the structure 
constant $C_{abc,\rho}^{\mu\nu}\in\C$ are obtained by direct calculations. 
\qed\end{pf}
In the same way as in the commutative ($\theta=0$) case, 
the product formula above leads to the following. 
Let $\Ob:=\{a,b,\cdots\}$ be a finite collection of labels, 
where any $a\in\Ob$ is associated with a nondegenerate 
symmetric matrix $A_a\in\Mat_n(\Z)$ such that 
for any $a,b\in\Ob$ the condition (\ref{DGcd}) holds and 
$A_{ab}$ is nondegenerate if $a\ne b$. 
For any $a,b\in\Ob$, define a vector space $H^0(a,b)$ as
follows: 
\begin{itemize}
 \item 
If $A_{ab}$ is positive definite, $H^0(a,b)$ is the 
$|\det(A_{ab})|$-dimensional vector space spanned by the theta functions
$\{\fe_{ab}^\mu\}$. 

 \item If $a=b$, then $H^0(a,b):=\C$. 

 \item If otherwise, then we set $H^0(a,b)=0$. 
\end{itemize}
Then, the product formula in Theorem \ref{thm:main} defines 
an algebraic structure on $\oplus_{a,b\in\Ob}H^0(a,b)$. 
The condition (\ref{DGcd}) has an interpretation in a categorical 
setting of these structures (see \cite{ncDG}).

 \section{An example}
\label{sec:exmp}

We end with showing an example for the case of noncommutative 
complex two-torus ($n=2$). 
In this case, for any fixed $\theta$, 
the condition $A_a\theta A_a=A_b\theta A_b$ reduces to 
\begin{equation*}
 \det(A_a)=\det(A_b)\ .
\end{equation*}
In general there exist infinite 
symmetric matrices $A\in\Mat_2(\Z)$ for a fixed $\det(A)$. 
For instance, diagonal matrices $A\in\Mat_2(\Z)$ with $\det(A)=-4$ are 
\begin{equation*}
 A_{\bf 1}=\bp 1 & 0 \\ 0 & -4\ep\ ,\quad 
 A_{\bf 2}=\bp 2 & 0 \\ 0 & -2\ep\ ,\quad 
 A_{\bf 3}=\bp 4 & 0 \\ 0 & -1\ep\ ,
\end{equation*}
and $A_{\bf 1'}:=-A_{\bf 1}$, $A_{\bf 2'}:=-A_{\bf 2}$, 
$A_{\bf 3'}:=-A_{\bf 3}$. 
Since $H^0(i,j')=H^0(i',j)=0$ for any $i,j={\bf 1},{\bf 2},{\bf 3}$, 
let us concentrate on the one side $\{{\bf 1},{\bf 2},{\bf 3}\}$. 
Then, one obtains $H^0(i,j)\ne 0$ if and only if $i\le j$ and 
in particular 
\begin{equation*}
 \dim(H^0({\bf 1},{\bf 2}))=2\ ,\qquad 
 \dim(H^0({\bf 2},{\bf 3}))=2\ ,\qquad 
 \dim(H^0({\bf 1},{\bf 3}))=9\ .
\end{equation*}
Thus, one obtains the following quiver: 
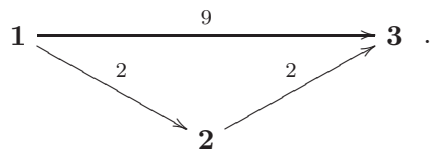
\begin{equation*}
\xymatrix{
 {\bf 1} \ar[drr]^{2} \ar[rrrr]^{9} &&&& {\bf 3} \\
 &&  {\bf 2} \ar[urr]^{2} && 
}\ .
\end{equation*}
However, there exist infinite 
symmetric matrices $A\in\Mat_n(\Z)$ with $\det(A)=-4$, 
since the matrix $g^tA g$ has $\det(A)=-4$ for any $SL(2,\Z)$ element $g$.


\end{document}